\documentclass[12pt,twoside]{article}
\usepackage{amsmath,amsbsy,amsfonts}

\newtheorem{theorem}{Theorem}[section]

\let\Section=\section
\def\section{\setcounter{equation}{0}\Section}
\begin{document}

\date{}
\title{Existence of solution for a class of nonlocal elliptic problem
via sub-supersolution method}
\author{\textsf{{Claudianor O. Alves} \thanks{%
C.O. Alves was partially supported by CNPq/Brazil 303080/2009-4,
coalves@dme.ufcg.edu.br}} \\
{\small \textit{Unidade Acad\^emica de Matem\'atica }}\\
{\small \textit{Universidade Federal de Campina Grande}}\\
{\small \textit{58429-900, Campina Grande - PB - Brazil}}\\
{\small \textit{e-mail address: coalves@dme.ufcg.edu.br}} \\
\\
\vspace{1mm} \textsf{{Dragos-Patru Covei} }\\
{\small \textit{Department of Applied Mathematics}}\\
{\small \textit{The Bucharest University of Economics Study }}\\
{\small \textit{Piata Romana, 1st district, postal code: 010374, postal
office: 22, Romania}}\\
{\small \textit{e-mail address:coveidragos@yahoo.com}}}
\maketitle

\begin{abstract}
We show the existence of solution for some classes of nonlocal problems. Our
proof combines the presence of sub and supersolution with the pseudomonotone operators theory.
\end{abstract}



\section{Introduction}

In this paper, we establish existence of solution for the
following class of nonlocal elliptic problem 
\begin{equation}
\left\{ 
\begin{array}{l}
 -a(\int_{\Omega }|u|^{q})\Delta {u}=h_1(x,u)f(\int_{\Omega
}|u|^{p})+h_2(x,u)g(\int_{\Omega}|u|^{r}),\,\,\,\mbox{in}\,\,\Omega, \\ 
\mbox{} \\ 
u=0 ,\,\,\,\mbox{on}\,\,\,\partial \Omega,%
\end{array}%
\right.  \label{1}
\end{equation}%
where $\Omega \subset \mathbb{R}^{N}\left( N\geq 1\right) $ is a smooth bounded domain, $q,p,r \in \lbrack 1,+\infty
) $, $h_i:\overline{\Omega }\times \mathbb{R}^{+} \rightarrow \mathbb{R}$ and $a$, $f$, $g:[0,+\infty )\rightarrow
(0,+\infty )$ are continuous functions with $f$, $g\in L^{\infty
}([0,+\infty ))$ and 
\begin{equation} 
\inf_{t\in \lbrack 0,+\infty )}a(t),\inf_{t\in \lbrack 0,+\infty
)}f(t),\inf_{t\in \lbrack 0,+\infty )}g(t)\geq a_{0}>0.  \label{(3.1)}
\end{equation}

The interest of such problems come from the articles of Chipot \& Lovat \cite%
{CL,CL2}, Chipot \& Rodrigues \cite{CRO}, Chipot \& Corr\^ea \cite{CJ} and Corr\^{e}a, Menezes \&
Ferreira \cite{CSJ}, where the authors study classes of nonlocal problems
motivated by the fact that they appear in some applied mathematics areas.
More exactly, it is pointed out in the paper \cite[see pp. 4619-4620 ]{CL}
that if $q=1$ and $h$ is of the form $h(x,u)$, the solution $u$ of the
problem (\ref{1}) could describe the density of a population
subject to spreading where the diffusion coefficient $a$ is supposed to
depend on the entire population in the domain rather than on the local
density. Moreover, in \cite{CL}, the authors have mentioned that the importance of such model lies in the fact 
that measurements that serves to determine physical constants are not made at a 
point but represent an average in a neighborhood of a point so that
these physical constants depend on local averages.  

It is worthwhile to remind that there are other classes of problems in which
nonlocal terms appear. For instance, in mechanics arises equations like
$$
\left\{
\begin{array}{l}
-M(\|u\|^{2})\Delta{u} = f(u), \quad \mbox{in} \quad \Omega, \\
\mbox{}\\
u=0, \quad \mbox{on} \quad \partial \Omega,
\end{array}
\right.
$$
related with nonlinear vibrations of beams, where $M:\mathbb{R} \to \mathbb{R}$ is a given function and $\|\,\,\,\|$ denotes the usual norm in $H_{0}^{1}(\Omega)$. More details about this class of problem can be see in \cite{AC2}, \cite{AC3}, \cite{CF}, \cite{F}, \cite{L} and \cite{Ma} and their references.

Our main interest in this work is the studying of the nonlocal problem (\ref{1}%
) for nonlinearities and boundary conditions which were not considered yet
in the literature. For example, in Section 3, we will consider the existence of solution for the class
of singular nonlocal problem 
\begin{equation*}
\left\{ 
\begin{array}{l}
-a(\int_{\Omega }|u|^{q})\Delta {u}=f(\int_{\Omega }|u|^{p}){u^{-\alpha }}%
+g(\int_{\Omega }|u|^{r})u^{\beta },\,\,\,\mbox{in}\,\,\Omega, \\ 
\mbox{} \\ 
u(x)>0\,\,\mbox{in}\,\,\Omega, \\ 
u=0\,\,\,\mbox{on}\,\,\,\partial \Omega,%
\end{array}%
\right. \eqno(SP)
\end{equation*}%
with $\alpha \in (0,1)$ and $\beta \in \lbrack 1,+\infty )$.

In Section 4,  we will prove the
existence of solution for the following class of nonlocal problem with
positive power 
\begin{equation*}
\left\{ 
\begin{array}{l}
-a(\int_{\Omega }|u|^{q})\Delta {u}=f(\int_{\Omega }|u|^{p}dx)u^{\alpha
}-g(\int_{\Omega }|u|^{r}dx)u^{\beta },\,\,\,\mbox{in}\,\,\Omega , \\ 
\mbox{} \\ 
u(x)>0,\,\,\mbox{in}\,\,\Omega ,\\ 
u=0,\,\,\,\mbox{on}\,\,\,\partial \Omega, %
\end{array}%
\right. \eqno(DP)
\end{equation*}%
with $\alpha \in [0,1)$ and $\beta \in \lbrack 1,+\infty )$.

The main tool used in the study of the above problems is the
sub-supersolution method. However, since we do not assume any hypotheses
involving monotonicity of the functions $a$, $f$ and $g$, we cannot used the
sub-supersolution method combined with maximum principle. Here, motivated by some arguments  found in Leon \cite{leon}, we prove a
new result combining the existence of sub and supersolution with the pseudomonotone operators theory, which can be used to study a large class of nonlocal problem. Our main result is the following 

\begin{theorem}
\label{T1} Let $\Omega \subset \mathbb{R}^{N}\left( N\geq 1\right) $ be a
smooth bounded domain and $p,q,r \in \lbrack 1,+\infty )$. Suppose that $h_i:\overline{\Omega }\times \mathbb{R} \rightarrow \mathbb{R}^{+}$ and $a$, $f$, $g:[0,+\infty )\rightarrow (0,+\infty )$ are continuous functions verifying (\ref{(3.1)}) with $f$, $g\in L^{\infty
}([0,+\infty ))$. Assume also that there are $\underline{u}$, $\overline{u}\in H^{1}(\Omega )\cap L^{\infty }(\Omega )$ satisfying the following inequalities:
\begin{equation} 
\underline{u}(x)\leq 0 \leq \overline{u}(x)\,\,\,\mbox{on}%
\,\,\,\partial \Omega,
\end{equation}
\begin{equation}
- \Delta \overline{u} \geq  \frac{1}{a_0}(h_1(x,\overline{u})\|f\|_{\infty}+ h_2(x,\overline{u})\|g\|_{\infty}), \quad \mbox{in} \quad \Omega  \label{s2}
\end{equation}
and 
\begin{equation}
-\Delta \underline{u} \leq  \frac{a_0}{\Gamma}(h_1(x,\underline{u})+ h_2(x,\underline{u})), \quad \mbox{in} \quad \Omega \label{s1}
\end{equation}
where
$$
\Gamma=\max\{a(t)\,:\, t \in [0, \|\overline{u}\|_{\infty}^{q}|\Omega|]\}.
$$
Then, there exists $u\in H_{0}^{1}(\Omega )\cap L^{\infty }(\Omega )$ 
with $\underline{u}(x)\leq u(x)\leq \overline{u}(x)\,\,\,\forall x\in \Omega
$ , which verifies  
$$
\begin{array}{l}
 a(\int_{\Omega}|u|^{q}){\displaystyle \int_{\Omega }}\nabla u\,\nabla \phi  =f(\int_{\Omega}|{u}|^{p}){\displaystyle \int_{\Omega
}}h_1(x,{u})\phi + g(\int_{\Omega}|{u}|^{r})  {\displaystyle \int_{\Omega}}h_2(x,{u})\phi \, \\
\end{array}
$$
for all $\phi \in H_{0}^{1}(\Omega )$, that is, $u$ is a weak solution of
the nonlocal problem 
$$
\left\{ 
\begin{array}{l}
- a(\int_{\Omega}|u|^{q}) \Delta {u}=h_1(x,{u})f(\int_{\Omega}|{u}|^{p})+ h_2(x,{u})g(\int_{\Omega}|{u}|^{r}),\,\,\,\mbox{in}%
\,\,\,\Omega , \\ 
\mbox{} \\ 
u(x)=0,\,\,\,\mbox{on}\,\,\,\partial \Omega .%
\end{array}%
\right. \eqno{(NLP)}
$$
\end{theorem}

\section{Proof of Theorem \protect\ref{T1}}

In what follows, we will consider the following functions
$$
\zeta(x,t)=
\left\{
\begin{array}{l}
\underline{u}(x), \quad \mbox{if} \quad t \leq \underline{u}(x)\\
t, \quad \mbox{if} \quad \underline{u}(x) \leq t \leq \overline{u}(x)\\
\overline{u}(x), \quad \mbox{if} \quad t \geq \overline{u},
\end{array}
\right.
$$ 
$$
\tilde{h}_i(x,t)=
\left\{
\begin{array}{l}
h_i(x,\underline{u}(x)), \quad \mbox{if} \quad t \leq \underline{u}(x)\\
h_i(x,t), \quad \mbox{if} \quad \underline{u}(x) \leq t \leq \overline{u}(x)\\
h_i(x,\overline{u}(x)), \quad \mbox{if} \quad t \geq \overline{u},
\end{array}
\right.
$$ 
and for $l \in (0,1)$, the function 
$$
\gamma(x,t)=-(\underline{u}(x)-t)_{+}^{l}+(t-\overline{u}(x))_{+}^{l}.
$$
Using the above functions, we study the existence of solution for the ensuing auxiliary nonlocal  problem 
\begin{equation*}
\left\{ 
\begin{array}{l}
-a(\int_{\Omega}|\zeta(x,u)|^{q})\Delta {u}=H(x,u,f(\int_{\Omega }|\zeta(x,u)|^{p}),g(\int_{\Omega }|\zeta(x,u)|^{r})),\,\,\,%
\mbox{in}\,\,\,\Omega , \\ 
\mbox{} \\ 
u=0 ,\,\,\,\mbox{on}\,\,\,\partial \Omega .%
\end{array}%
\right. \eqno{(AP)}
\end{equation*}%
where
$$
H(x,u,s,t)=\tilde{h}_1(x,u)s+\tilde{h}_2(x,u)t- \gamma(x,u).
$$
To prove the existence of solution for problem $(AP)$, we will use the pseudomonotone operators theory. To this end, we will work with the operator $B:H^{1}_{0}(\Omega) \to H^{-1}$ given by
$$
\begin{array}{l}
\left\langle B(u),v\right\rangle=  a(\int_{\Omega}|\zeta(x,u)|^{q}) \displaystyle  \int_{\Omega} \nabla u \nabla v \\
\mbox{}\\
\;\;\;\;\;\;\;\;\;\;\;\;\;\;\;\;\;\;- {\displaystyle \int_{\Omega}}H(x,u,f(\int_{\Omega }|\zeta(x,u)|^{p}),g(\int_{\Omega }|\zeta(x,u)|^{r}))v \quad \forall u,v \in H^{1}_{0}(\Omega).
\end{array}
$$ 
An direct computation gives that $B$ is continuous, bounded and coercive, i.e., 
$$
\frac{\left\langle B(u),u \right\rangle}{\|u\|} \to +\infty \quad \mbox{as} \quad \|u\| \to +\infty.
$$
Moreover, $B$ is a pseudomonotone operator, i.e., if $(u_n) \subset H^{1}_{0}(\Omega)$ verifies $u_n \rightharpoonup u$ in $H^{1}_{0}(\Omega)$ and 
\begin{equation} \label{B1}
\limsup_{n \to +\infty} \left\langle B(u_n), u_n-u \right\rangle \leq 0,
\end{equation}
then 
\begin{equation} \label{B22}
\liminf_{n \to +\infty} \left\langle B(u_n), u_n-v \right\rangle \geq \left\langle B(u), u-v \right\rangle   \quad \forall v \in H^{1}_{0}(\Omega).
\end{equation}
In fact, by using the definition of $H$ combined with Lebesgue's Theorem, it follows that
$$
\frac{{\displaystyle \int_{\Omega}}H(x,u_n,f(\int_{\Omega }|\zeta(x,u_n)|^{p}),g(\int_{\Omega }|\zeta(x,u_n)|^{r}))(u_n-u)}{a(\int_{\Omega}|\zeta(x,u_n)|^{q})} \to 0.
$$
The last limit together with (\ref{B1}) leads to
\begin{equation} \label{B3}
\limsup_{n \to +\infty}\int_{\Omega}\nabla u_n \nabla (u_n-u) \leq 0.
\end{equation}
Once that
$$
\int_{\Omega}\nabla u \nabla (u_n-u) \to 0
$$
and
$$
\|u_n-u\|^{2}=\int_{\Omega}\nabla u_n \nabla(u_n-u)-\int_{\Omega}\nabla u \nabla(u_n-u),
$$
we can conclude that
$$
\|u_n-u\|^{2} \to 0,
$$
or equivalently, 
\begin{equation} \label{B33}
u_n \to u \quad \mbox{in} \quad  H^{1}_{0}(\Omega).
\end{equation}
From this, we observe that (\ref{B22}) is an immediate consequence of (\ref{B33}).

From \cite[Theorem 3.3.6]{Necas}, $B$ is surjective, i.e., $B(H^{1}_{0}(\Omega))=H^{-1}$. Therefore,  there exists $u \in H^{1}_{0}(\Omega)$ such that 
$$
\left\langle B(u),\phi \right\rangle=0 \,\,\ \forall \phi \in E,
$$ 
implying that $u$ is a solution of $(AP)$. Now, our goal is showing that $u$ is a solution of $(NLP)$. To this end, we must prove that
\begin{equation} \label{z1}
\underline{u} \leq u \leq \overline{u} \quad \mbox{in} \quad \Omega. 
\end{equation}
Choosing $\phi=(u-\overline{u})_{+}$ as test function, we derive
$$
a({\textstyle \int_{\Omega}}|u|^{q}) \int_{\Omega}\nabla u \nabla (u-\overline{u})_{+}=\int_{\Omega}H(x,u,f( {\textstyle \int_{\Omega }}|\zeta(x,u)|^{p})),g({\textstyle \int_{\Omega }}|\zeta(x,u)|^{r})))(u-\overline{u})_{+}
$$
that is,
$$
\begin{array}{l}
a(\int_{\Omega}|u|^{q}){{\displaystyle \int_{\Omega}}}\nabla u \nabla (u-\overline{u})_{+}=f( \int_{\Omega }|\zeta(x,u)|^{p})){\displaystyle \int_{\Omega}}h_1(x,u)(u-\overline{u})_{+} \\
\mbox{}\\
\hspace{4 cm} + g( \int_{\Omega }|\zeta(x,u)|^{r})){\displaystyle \int_{\Omega}}h_2(x,u)(u-\overline{u})_{+} - {\displaystyle \int_{\Omega}}\gamma(x,u)(u-\overline{u})_{+}.
\end{array}
$$
Using the definition of $h_i\,(i=1,2)$ and $\gamma$, it follows that
$$
\begin{array}{l}
a(\int_{\Omega}|u|^{q}){\displaystyle \int_{\Omega}}\nabla u \nabla (u-\overline{u})_{+}=f( \int_{\Omega }|\zeta(x,u)|^{p})){\displaystyle \int_{\Omega}}h_1(x,\overline{u})(u-\overline{u})_{+} \\
\mbox{}\\
\hspace{5 cm} +g( \int_{\Omega }|\zeta(x,u)|^{r})){\displaystyle \int_{\Omega}}h_2(x,\overline{u})(u-\overline{u})_{+} - {\displaystyle \int_{\Omega}}(u-\overline{u})_{+}^{l+1}.
\end{array}
$$
Thus,
$$
\begin{array}{l}
{\displaystyle \int_{\Omega}}\nabla u \nabla (u-\overline{u})_{+} \leq {\displaystyle \frac{1}{a_0}}\|f\|_{\infty}{\displaystyle \int_{\Omega}}h_1(x,\overline{u})(u-\overline{u})_{+} \\
\mbox{}\\
\hspace{3 cm} +{\displaystyle \frac{1}{a_0}}\|g\|_{\infty}{\displaystyle \int_{\Omega}}h_2(x,\overline{u})(u-\overline{u})_{+} - {\displaystyle \frac{1}{a(\int_{\Omega}|u|^{q})}}{\displaystyle \int_{\Omega}}(u-\overline{u})_{+}^{l+1}.
\end{array}
$$
From (\ref{s1}), 
$$
\int_{\Omega} \nabla u \nabla (u-\overline{u})_{+} \leq \int_{\Omega} \nabla \overline{u} \nabla (u-\overline{u})_{+} - \frac{1}{a(\int_{\Omega}|u|^{q})}{\displaystyle \int_{\Omega}}(u-\overline{u})_{+}^{l+1},
$$
or equivalently,
$$
\|(u- \overline{u})_{+}\|^{2} \leq - \frac{1}{a(\int_{\Omega}|u|^{q})}{\displaystyle \int_{\Omega}}(u-\overline{u})_{+}^{l+1} \leq 0,
$$
showing that $(u-\overline{u})_+=0$.

Now, to prove that $\underline{u} \leq u$, we choose $\phi=(\underline{u}-u)_{+}$ as test function. Repeating the above arguments, we obtain 
$$
\begin{array}{l}
a(\int_{\Omega}|u|^{q}){\displaystyle \int_{\Omega}}\nabla u \nabla (\underline{u}-u)_{+}=f( \int_{\Omega }|\zeta(x,u)|^{p})){\displaystyle \int_{\Omega}}h_1(x,\overline{u})(\underline{u}-u)_{+} \\
\mbox{}\\
\hspace{5 cm} g( \int_{\Omega }|\zeta(x,u)|^{r})){\displaystyle \int_{\Omega}}h_2(x,\overline{u})(\underline{u}-u)_{+} + {\displaystyle \int_{\Omega}}(\underline{u}-u)_{+}^{l+1}.
\end{array}
$$
Thus, 
$$
\begin{array}{l}
{\displaystyle \int_{\Omega}}\nabla u \nabla (\underline{u}-u)_{+} \geq  \displaystyle \frac{a_0}{\Gamma}{\displaystyle \int_{\Omega}}h_1(x,\overline{u})(\underline{u}-u)_{+} \\
\mbox{}\\
\hspace{5 cm} \displaystyle \frac{a_0}{\Gamma}{\displaystyle \int_{\Omega}}h_2(x,\overline{u})(\underline{u}-u)_{+} + {\displaystyle \frac{1}{\Gamma}\int_{\Omega}}(\underline{u}-u)_{+}^{l+1}.
\end{array}
$$
From (\ref{s2}), 
$$
\int_{\Omega} \nabla u \nabla (\underline{u}-u)_{+} \geq \int_{\Omega} \nabla \underline{u} \nabla (\underline{u}-u)_{+} + \frac{1}{\Gamma}{\displaystyle \int_{\Omega}}(\underline{u}-u)_{+}^{l+1},
$$
or equivalently,
$$
\|(\underline{u}-u)_{+}\|^{2} \leq - \frac{1}{\Gamma}{\displaystyle \int_{\Omega}}(\underline{u}-u)_{+}^{l+1} \leq 0,
$$
showing that $(\underline{u}-{u})_+=0$. Therefore, Theorem \ref{T1} is proved. \hfill \rule{2mm}{2mm}

\section{Application I : Existence of solution for a class of nonlocal
problem with singular term\label{4}}

In this section, we establish the existence of solution for the
following class of singular nonlocal problem 
\begin{equation*}
\left\{ 
\begin{array}{l}
- a(\int_{\Omega }|u|^{q})\Delta {u}=\frac{f( \int_{\Omega
}|u|^{p})}{u^{\alpha }}+g(\int_{\Omega }|u|^{r})u^{\beta },\,\,\,%
\mbox{in}\,\,\Omega, \\ 
\mbox{} \\ 
u(x)>0,\,\,\mbox{in}\,\,\Omega, \\ 
u=0,\,\,\,\mbox{on}\,\,\,\partial \Omega,%
\end{array}%
\right. \eqno(SP)
\end{equation*}%
where $\Omega $ is a bounded domain with smooth boundary and $q$, $p$, $r\in
\lbrack 1,+\infty )$, $\alpha \in \lbrack 0,1)$, $\beta \geq 0$ and $a$, $f$%
, $g:[0,+\infty )\rightarrow \mathbb{R}^{+}$ are continuous functions verifying
condition (\ref{(3.1)}) with $f$, $g\in L^{\infty }([0,+\infty))$.

The existence of sub and supersolution verifying the hypotheses of Theorem \ref{T1} can be obtained of the following way: Fix $R>0$ such that $\Omega \subset B_{R}(0)$. If $e
\in C^{2}(\overline{\Omega})$ denotes the unique positive solution of 
\begin{equation*}
\left\{ 
\begin{array}{l}
-\Delta e = 1, \,\,\, \mbox{in} \,\,\ B_{R}(0), \\ 
\mbox{} \\ 
e=0, \,\,\, \mbox{on} \,\,\, \partial B_{R}(0),
\end{array}
\right.
\end{equation*}
and $M>0$ is a constant large enough, we derive 
\begin{equation*}
\left\{ 
\begin{array}{l}
-\Delta (Me) = M  \geq {\displaystyle \frac{1}{a_0}} (
\|f\|_{\infty}(Me)^{\alpha}+\|g\|_{\infty}(Me)^{\beta}), \,\,\, \mbox{in} \,\,\ \Omega ,\\ 
\mbox{} \\ 
M e >0, \,\,\, \mbox{on} \,\,\, \partial \Omega.
\end{array}
\right.
\end{equation*}
On the other hand, if $\phi_1$ denotes a positive eigenfunction associated with the first
eigenvalue $\lambda_1$ of $(-\Delta, H_{0}^{1}(\Omega))$, we observe that
for $\epsilon >0$ smaller enough, the below inequality occurs 
\begin{equation*}
-\Delta (\epsilon \phi_1) \leq \frac{a_0}{\Gamma}
\Big(\frac{1}{(|\epsilon \phi_1|^{2}+\delta )^{\frac{\alpha}{2}}}
+(\epsilon \phi_1)^{\beta}\Big) \,\,\, \mbox{in} \,\,\,
\Omega,
\end{equation*}
where, 
$$
\Gamma=\max\{a(t)\,:\, t \in [0, \|Me\|_{\infty}^{q}|\Omega|]\},
$$
for $\epsilon$ small enough uniformly for $ \delta \in (0,1]$.

Consequently, for $\epsilon$ small enough and $M$ large enough, the functions $\underline{u}=\epsilon \phi$ and $\overline{u}%
=Me$ verify the hypotheses of Theorem \ref{T1} for the nonlocal problem 
\begin{equation*}
\left\{ 
\begin{array}{l}
- a(\int_{\Omega }|u|^{q})\Delta {u}=\frac{f(\int_{%
\Omega}|u|^{p}dx)}{(|u|^{2}+\delta)^{\frac{\alpha}{2}}}+g(\int_{%
\Omega}|u|^{r}dx)u^{\beta},\,\,\,\mbox{in}\,\,\Omega ,\\ 
\mbox{} \\ 
u=0 \,\,\,\mbox{on}\,\,\,\partial \Omega .%
\end{array}%
\right. \eqno{(P_\delta)}
\end{equation*}
Thereby, there exists a solution $u_\delta \in H_{0}^{1}(\Omega)$ for the nonlocal 
problem 
\begin{equation*}
\left\{ 
\begin{array}{l}
-a(\int_{\Omega }|u_\delta|^{q})\Delta {u_\delta}=\frac{%
f(\int_{\Omega}|u_\delta|^{p}dx)}{(|u_\delta|^{2}+\delta)^{\frac{\alpha}{2}}}%
+g(\int_{\Omega}|u_\delta|^{r}dx)u_\delta^{\beta},\,\,\,\mbox{in}\,\,\Omega,
\\ 
\mbox{} \\ 
u_\delta(x) >0, \,\, \mbox{in} \,\, \Omega, \\ 
u_\delta=0, \,\,\,\mbox{on}\,\,\,\partial \Omega,%
\end{array}%
\right. \eqno{(P_\delta)}
\end{equation*}
with 
\begin{equation*}
\underline{u} \leq u_\delta \leq \overline{u} \,\,\, \mbox{in} \,\,\, \Omega
\,\,\, \forall \delta \in [0,1].
\end{equation*}
In what follows, for each $n \in \mathbb{N}$, we denote by $u_n$ the
solution $u_{\frac{1}{n}}$. Therefore 
\begin{equation*}
\left\{ 
\begin{array}{l}
-a(\int_{\Omega }|u_n|^{q})\Delta {u_n}=\frac{%
f(\int_{\Omega}|u_n|^{p})}{(|u_n|^{2}+\frac{1}{n})^{\frac{\alpha}{2}}}%
+g(\int_{\Omega}|u_n|^{r})u_n^{\beta},\,\,\,\mbox{in}\,\,\Omega, \\ 
\mbox{} \\ 
u_n(x) >0, \,\, \mbox{in} \,\, \Omega, \\ 
u_n=0 , \,\,\,\mbox{on}\,\,\,\partial \Omega,%
\end{array}%
\right. \eqno{(P_n)}
\end{equation*}
and 
\begin{equation}  \label{D2}
\underline{u} \leq u_n \leq \overline{u} \,\,\, \mbox{in} \,\,\, \Omega
\,\,\, \forall n \in \mathbb{N}.
\end{equation}
Once that $u_n$ is a solution of $(P_n)$, we have the below equality 
\begin{equation}  \label{B2}
a({\textstyle \int_{\Omega }}|u_n|^{q}dx)\int_{\Omega}\nabla u_n \nabla v dx =
\int_{\Omega}\frac{f(\int_{\Omega}|u_n|^{p})v}{(|u_n|^{2}+\frac{1}{n})^{%
\frac{\alpha}{2}}}dx+\int_{\Omega}g({\textstyle \int_{\Omega}}|u_n|^{r}dx)u_n^{\beta}v dx
\end{equation}
for all $v \in H^{1}_{0}(\Omega)$. Recalling that 
\begin{equation*}
a(t) \geq a_o \,\,\, \forall t \geq 0,
\end{equation*}
and $f,g \in L^{\infty}([0,+\infty))$, it follows that 
\begin{equation*}
a_0\int_{\Omega}|\nabla u_n|^{2} \leq C
\int_{\Omega}(u_n^{1-\alpha}+u_n^{\beta+1}) dx.
\end{equation*}
Now, using that $\alpha \in [0,1)$, $\beta \geq 0$ and (\ref{D2}), the last
inequality gives that $(u_n)$ is bounded in $H^{1}_{0}(\Omega)$. Thus, for
some subsequence, still denote by $(u_n)$, there exists $u \in
H_{0}^{1}(\Omega)$ such that 
\begin{equation*}
u_n \rightharpoonup u \,\,\, \mbox{in} \,\,\, H_{0}^{1}(\Omega)
\end{equation*}
and 
\begin{equation*}
u_n(x) \to u(x) \,\,\ \mbox{a.e in} \,\,\, \Omega.
\end{equation*}
Since 
\begin{equation*}
\underline{u} \leq u_n \leq \overline{u} \,\,\, \mbox{in} \,\, \Omega,
\end{equation*}
the last limit yields 
\begin{equation*}
u_n \to u \,\,\, \mbox{in} \,\,\, L^{s}(\Omega) \,\,\, \forall s \in
[1,+\infty).
\end{equation*}
Thus, by continuity of $a,f$ and $g$, 
\begin{equation}  \label{B3}
a({\textstyle \int_{\Omega}}|u_n|^{q}dx) \to a({\textstyle \int_{\Omega}}|u|^{q}dx),
\end{equation}
\begin{equation}  \label{B31}
f({\textstyle \int_{\Omega}}|u_n|^{p}dx) \to f({\textstyle \int_{\Omega}}|u|^{p}dx)
\end{equation}
and 
\begin{equation}  \label{B32}
g({\textstyle \int_{\Omega}}|u_n|^{r}dx) \to g({\textstyle \int_{\Omega}}|u|^{r}dx).
\end{equation}

Taking the limit in (\ref{B2}) with $v\in C_{0}^{\infty }(\Omega )$ and
using (\ref{B3})-(\ref{B32}), we get 
\begin{equation*}
a({\textstyle \int_{\Omega} }|u|^{q}dx)\int_{\Omega }\nabla u\nabla vdx=\int_{\Omega }%
\frac{f(\int_{\Omega }|u|^{p}dx)v}{u^{\alpha }}dx+\int_{\Omega
}g({\textstyle \int_{\Omega} }|u|^{r}dx)vdx.
\end{equation*}%
Now, repeating the same arguments explored in Alves \& Corr\^{e}a \cite[Page
735]{AC}, we can conclude that for all $v\in H_{0}^{1}(\Omega)$ 
\begin{equation}
a({\textstyle \int_{\Omega} }|u|^{q}dx)\int_{\Omega }\nabla u\nabla vdx=\int_{\Omega }%
\frac{f(\int_{\Omega }|u|^{p}dx)v}{u^{\alpha }}dx+\int_{\Omega
}g({\textstyle  \int_{\Omega }}|u|^{r}dx)vdx,  \label{B4}
\end{equation}%
showing that $u$ is a solution of nonlocal problem $(SP)$. From the above
commentaries, we have proved the following result

\begin{theorem}
\label{T3} Assume that $a,f,g:[0,+\infty) \to \mathbb{R}$ are continuous
functions verifying condition (\ref{(3.1)}) with $f,g \in
L^{\infty}([0,+\infty))$, $p,q,r \in [1,+\infty)$, $\beta \geq 0$ and $%
\alpha \in [0,1)$. Then, problem $(SP)$ has a solution.
\end{theorem}

The Theorem \ref{T3} is related to the papers due to Coclite \& Palmieri 
\cite{CP} and Zhang \& Yu \cite{ZY}, in the sense that, in these papers the
authors considered the existence of solution for the local case, that is, $%
a=f=g=1$.

\section{Application II: Existence of solution for a class of nonlocal
problem with Dirichlet boundary conditions\label{3}}

In this section, we study the existence of positive solution for the
following class of nonlocal problem 
\begin{equation*}
\left\{ 
\begin{array}{l}
- a(\int_{\Omega }|u|^{q})\Delta {u}=f(\int_{\Omega
}|u|^{p}dx)u^{\alpha }-g(\int_{\Omega }|u|^{r}dx)u^{\beta },\,\,\,\mbox{in}%
\,\,\Omega \\ 
\mbox{} \\ 
u(x)>0\,\,\mbox{in}\,\,\Omega \\ 
u=0\,\,\,\mbox{on}\,\,\,\partial \Omega%
\end{array}%
\right. \eqno(DP)
\end{equation*}%
where $\Omega $ is a bounded domain with smooth boundary, $q$, $p$, $r$, $%
\beta \in \lbrack 1,+\infty )$ and $\alpha \in \lbrack 0,1)$. Related to the
functions $a$, $f$, $g:[0,+\infty )\rightarrow \mathbb{R}$, we assume they
verify (\ref{(3.1)}) and $f, g\in L^{\infty }([0,+\infty ))$.

We intend to use again Theorem \ref{T1}, however once that we are considering a negative signal between the terms in the right side of the problem, we must make an adjustment in that result.  Here, we will use the following version

\begin{theorem}
\label{T2} Let $\Omega \subset \mathbb{R}^{N}\left( N\geq 1\right) $ be a
smooth bounded domain and  $p,q,r \in \lbrack 1,+\infty )$. Suppose that $h_i:\overline{\Omega }\times \mathbb{R} \rightarrow \mathbb{R}^{+}$ and $a$, $f$, $g:[0,+\infty )\rightarrow (0,+\infty )$ are continuous functions with $a,f$, $g\in L^{\infty
}([0,+\infty ))$ and  verifying (\ref{(3.1)}). Assume also that there are $\underline{u}$, $\overline{u}\in H^{1}(\Omega )\cap L^{\infty }(\Omega )$ satisfying the following inequalities:
\begin{equation} 
\underline{u}(x)\leq 0 \leq \overline{u}(x)\,\,\,\mbox{on}%
\,\,\,\partial \Omega,
\end{equation}
\begin{equation}
- \Delta \overline{u} \geq  \frac{1}{a_0}h_1(x,\overline{u})\|f\|_{\infty}-\frac{a_0}{\|a\|_{\infty}}h_2(x,\overline{u}), \quad \mbox{in} \quad \Omega  \label{s22}
\end{equation}
and 
\begin{equation}
-\Delta \underline{u} \leq  \frac{a_0}{\|a\|_{\infty}}h_1(x,\underline{u})- \frac{\|g\|_{\infty}}{a_0}h_2(x,\underline{u}), \quad \mbox{in} \quad \Omega \label{s11}
\end{equation}
Then, there exists $u\in H_{0}^{1}(\Omega )\cap L^{\infty }(\Omega )$ 
with $\underline{u}(x)\leq u(x)\leq \overline{u}(x)\,\,\,\forall x\in \Omega
$  and $u \in H_{0}^{1}(\Omega )$, which verifies  
$$
\begin{array}{l}
 a(\int_{\Omega}|u|^{q}){\displaystyle \int_{\Omega }}\nabla u\,\nabla \phi  =f(\int_{\Omega}|{u}|^{p}){\displaystyle \int_{\Omega
}}h_1(x,{u})\phi - g(\int_{\Omega}|{u}|^{r})  {\displaystyle \int_{\Omega}}h_2(x,{u})\phi \, \\
\end{array}
$$
for all $\phi \in H_{0}^{1}(\Omega )$, that is, $u$ is a weak solution of
the nonlocal problem 
$$
\left\{ 
\begin{array}{l}
- a(\int_{\Omega}|u|^{q}) \Delta {u}=h_1(x,{u})f(\int_{\Omega}|{u}|^{p})- h_2(x,{u})g(\int_{\Omega}|{u}|^{r}),\,\,\,\mbox{in}%
\,\,\,\Omega \\ 
\mbox{} \\ 
u(x)=0,\,\,\,\mbox{on}\,\,\,\partial \Omega .%
\end{array}%
\right. \eqno{(NLP)}
$$
\end{theorem}

\vspace{0.5 cm}

\noindent {\bf Proof.} \quad The proof follows repeating the same type of arguments used in the proof of Theorem \ref{T1}.  \hfill \rule{2mm}{2mm}

\vspace{0.5 cm}

First of all, we observe that if $M>0$ is large enough, the function $\overline{u}=M$ verifies the below inequality 
\begin{equation*}
\left\{ 
\begin{array}{l}
-\Delta {\overline{u}}\geq  \frac{\|f\|_{\infty}}{a_0}\overline{u}^{\alpha}-\frac{a_0}{\|a\|_{\infty}}\overline{u}^{\beta},\,\,\,\mbox{in}\,\,\,\Omega \\ 
\mbox{} \\ 
\overline{u}>0,\,\,\,\mbox{on}\,\,\,\partial \Omega .%
\end{array}%
\right.
\end{equation*}

On the other hand, by a direct computation, if $\epsilon >0$ is smaller enough and $\phi_1$
denotes a positive eigenfunction associated with the first eigenvalue $%
\lambda_1$ of $(-\Delta , H^{1}_{0}(\Omega))$, it is easy to check that $%
\underline{u}=\epsilon \phi _{1}$ verifies 
\begin{equation*}
\left\{ 
\begin{array}{l}
-\Delta {\underline{u}}\leq
\frac{a_0}{\|a\|_{\infty}}\underline{u}^{\alpha}- \frac{\|g\|_{\infty}}{a_0}\underline{u}^{\beta}\,\,\,\,\mbox{in}\,\,\,\Omega \\ 
\mbox{} \\ 
\underline{u}=0,\,\,\,\mbox{on}\,\,\,\partial \Omega .%
\end{array}%
\right.
\end{equation*}%

From the above considerations, the functions $\underline{u}=\epsilon \phi $
and $\overline{u}=M$ verify the hypotheses of Theorem \ref{T1} for $\epsilon$ smaller enough and $M$ large enough. Thus, 
there exists $u\in H_{0}^{1}(\Omega )$ solution of $(DP)$. From the above commentaries, we have proved the ensuing result
\begin{theorem}
\label{T2} Assume that $a,f,g:[0,+\infty)\rightarrow \mathbb{R}$ are
continuous functions verifying condition (\ref{(3.1)}) with $f$, $g\in
L^{\infty }([0,+\infty))$, $p$, $q$, $r$, $\beta \in \lbrack 1,+\infty )$
and \linebreak $\alpha \in \lbrack 0,1)$. Then, problem $(DP)$ has a
solution.
\end{theorem}

The Theorem \ref{T2} is related with some results found in Alama \&
Tarantello \cite{AT}, Radulescu \& Repovs \cite{RR} and Lane \cite{Lane},
where the existence of solution for $(DP)$ have been considered
for the local case, that is, when $a, f, g=1$.

\end{document}